\documentclass[11pt]{article}

\usepackage{float}
\usepackage{amsmath}
\usepackage{amssymb}
\usepackage{array}
\usepackage[table]{xcolor} 
\usepackage{longtable}
\usepackage{fancyhdr}
\usepackage{pdflscape}
\usepackage{color}
\usepackage{colortbl}
\usepackage{rotating}
\usepackage{verbatim}
\usepackage{xcolor}
\usepackage[framemethod=TikZ]{mdframed}
\usepackage{lipsum}
\mdfdefinestyle{MyFrame}{
    linecolor=blue,
    outerlinewidth=1.0pt,
    roundcorner=20pt,
    innertopmargin=\baselineskip,
    innerbottommargin=\baselineskip,
    innerrightmargin=20pt,
    innerleftmargin=40pt,
    backgroundcolor=yellow!20!white}

\definecolor{green}{rgb}{0,1,0}
\definecolor{gray}{gray}{0.9}
\definecolor{LightCyan}{rgb}{0.88,1,1}
\definecolor{darkblue}{rgb}{0.0, 0.0, 0.5}
\definecolor{myLightGray}{gray}{0.6}

\fancyheadoffset{0.0in}

\begin{document}

\thispagestyle{empty}
\begin{center}
\Large{\color{blue}\bf A Robust Binary Nonlinear Solver \\ for Multi-stage Decisions } \\
\vspace{8mm}
\large{Mathematical Note} \\ 
\vspace{1cm}
\vspace{\baselineskip}
\small{Kashif Rashid}\footnote{\scriptsize{Formerly, Scientific Advisor, Schlumberger-Doll Research, One Hampshire Street, Cambridge, MA 02139. \\ Contact email: kash4questions@yahoo.com }} \\ 
\vspace{2mm}
\footnotesize{Scientist } \\
\vspace{8mm}
\small{May 2026} \\
\end{center}
\vspace{3cm}

\small
\section*{\sffamily \large\color{darkblue} Abstract }

This document describes two developments. The first is a Multi-stage Decision Framework that permits the formulation of various sequential decision-making problems in a generic manner. The second concerns a robust stochastic solution procedure designed to handle the binary nonlinear optimization problem stemming from the proposed decision framework. The method adopts Monte-Carlo search with incremental probabilistic learning with reinforcement. The two developments are described herein along with a demonstrative application example.

\vspace{3mm}

{\scriptsize\noindent {\emph{Key words:} Multi-stage, decisions, binary systems, probabilistic and reinforcement learning.  } }

\clearpage
\rhead{\tiny DECISION FRAMEWORK}
\small
\section*{\sffamily \large\color{darkblue} Multi-stage Decision Framework }

Multi-stage decision problems, with or without availability of new information, are not uncommon in real-world settings. These concern situations in which decisions are made at well-defined points (or stages), considered either spatially or temporally, to achieve a desired outcome over the entire horizon subject to local and global constraints~\cite{POWELL}. In this setting, local constraints apply to individual stages or components, while global constraints may apply to some metric over the entire decision horizon or some part thereof. The generic decision problem can be posed as one in which there are $N$ stages, each with $M$ components or choices. This gives rise to $M^N$ combinations or possible decision paths. Evidently the problem quickly becomes computationally involved with an increasing size for $N$ or $M$. 

\vspace{3mm}

The solution process serves to identify the optimal path (the selection of components in each stage) so as to maximize the desired objective. Consider the classic knapsack problem as an example. Here, given a collection of $N$ items, each with a given weight and utility, the goal is to select a group of items that maximize the total utility subject to a permitted weight limit~\cite{KNAP}. With only two choices per step ($M$=2), to select or not to select, this gives rise to $2^N$ possible outcomes. This number grows significantly with increasing $N$. While, efficient methods exist for mixed-integer linear problems of this type, for nonlinear cases in the general, the problem quickly becomes intractable due to its combinatorial nature. This is further exacerbated with an increasing number of choices $M$. Such problems may concern the selection of unit operations in a process train~\cite{SUPER}. Here the notion of the super-structure is employed in which all design choices are presented by stage and an optimal selection is sought to maximize yield, profit or other intended goal~\cite{YEO, MENCAR, KENKEL}. Mixed-integer nonlinear formulations are used in which component nonlinearity is invariably assuaged in preference for ease of resolution. Nonetheless, well-structured formulations are necessary to capture design optionality and restrictions using generalized disjunctive programming~\cite{ADJIMAN, GDP}. Such formulations are non-trivial and particularly cumbersome to develop for large-scale models comprising many design options. Thus, efficient solution methods are not only desirable, but necessary for highly nonlinear problems. This serves as the motivation for the multi-stage decision framework and the robust binary nonlinear solution approach.

\vspace{3mm}

\begin{figure}[ht]
    \centering
    \includegraphics[width=0.5\textwidth]{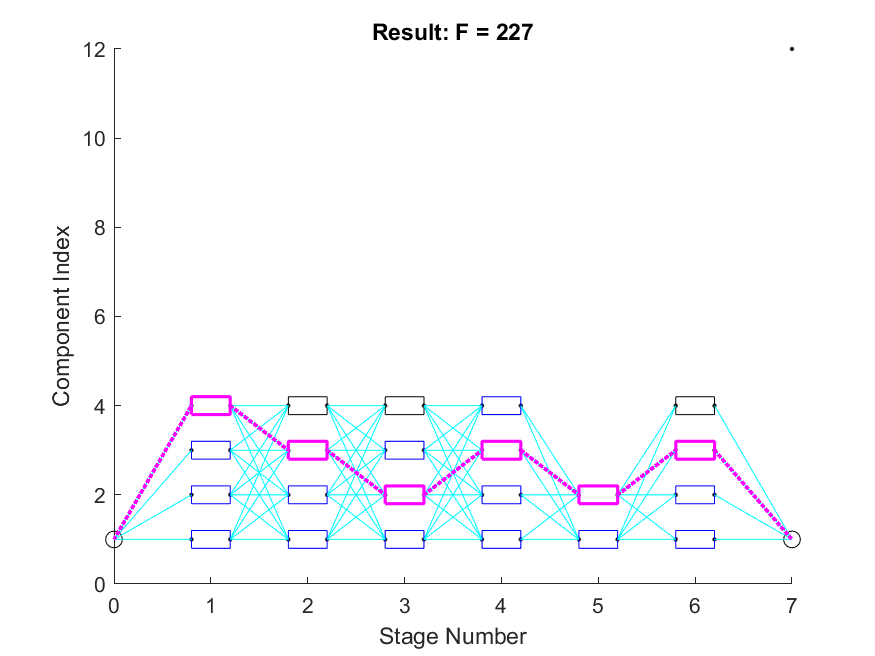}
    \caption{Solution path over 6 stages - fully connected case.} 
	\footnotesize{\label{fig:FullLink} 6 stages with [4 3 3 4 2 3] components per stage. The binary system is $2^{11}$ with 2048 choices.} \\
	\footnotesize{Blue squares represent valid choices, black squares indicate invalid choices and the edges show feasible links.}
    \label{fig:my_figure}
\end{figure}

\begin{figure}[ht]
    \centering
    \includegraphics[width=0.5\textwidth]{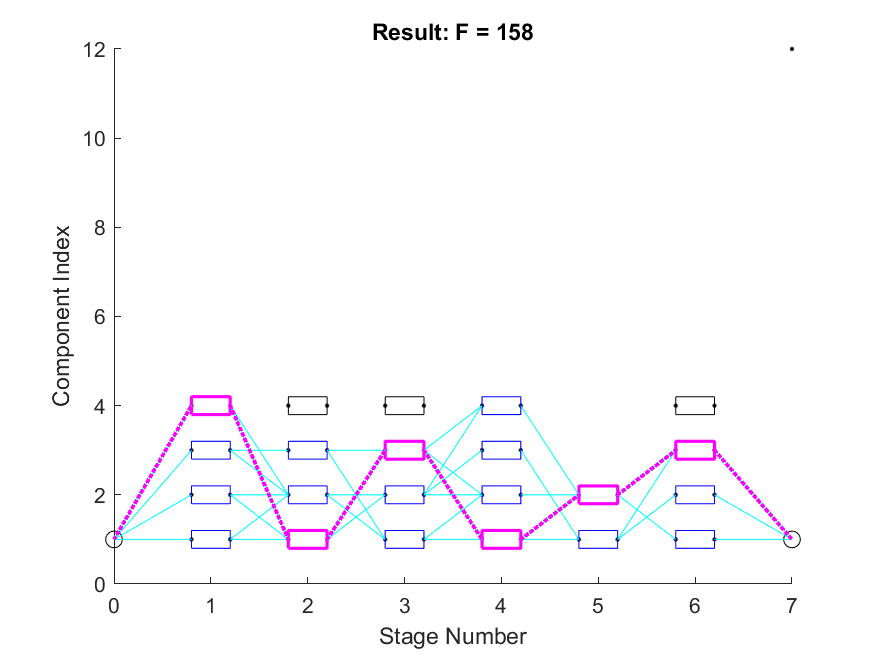}
    \caption{Solution path over 6 stages - partially connected case.}
	\footnotesize{\label{fig:PartialLink} 6 stages with [4 3 3 4 2 3] components per stage. The binary system is $2^{11}$ with 2048 choices (not all feasible).} \\
	\footnotesize{Blue squares represent valid choices, black squares indicate invalid choices and the edges show feasible links.}
    \label{fig:my_figure}
\end{figure}

In the multi-stage decision framework, the posed problem is reduced to a graph. A starting node connects to the nodes in the first layer, each of these will connect to those in the next layer and so on, until the nodes in the last layer connect to one terminating node (example shown in Figure \ref{fig:FullLink}). Thus, the goal is to find an optimal path from the start to the final node. The framework permits a fully connected graph (a node in one layer connects to all the nodes in the next layer) or a partially connected graph (where the feasible node connections are stipulated in advance) as per Figure \ref{fig:PartialLink}. In this regard, the graph may be initialized as fully connected and all undesired links are subsequently removed. Alternately, the graph may commence with zero connections and each permitted link is added in turn. We refer to a connection map in which node $i$ in layer ($j-1$) connects to node $k$ in layer $j$ if the flag is true, given by the tuple [$j-1$ $i$ $j$ $k$ 1], or indicates a deactivation if it is false, with [$j-1$ $i$ $j$ $k$ 0]. In this manner, any arrangement can be posed for consideration. For a sparse graph, only the feasible links need be described. Note that the scheme is flexible to permit the arrangement of connected sub-graphs or other configurations by need. The only stipulation is that the graph is directed (from left to right) and that the nodes connect layer-by-layer without jumps or return paths. Indeed, these elements, should they exist, can be removed from the outset by suitable graph construction. That is, with the insertion of dummy nodes and viable links as appropriate. Thus, many problems, indicative of restricted or conditioned choice, can be posed and formulated on the underlying directed graph.

\vspace{3mm}

Each node in a layer represents a selectable component. In the foregoing, this was noted as a simple connection in that the node had no other purpose than to indicate the selected path, where the objective value is a function of the chosen components. Now, each node can be described in further detail. We can assert the notion of a functional with integer input $y$ given by $v=f(y)$, where $y=1$ indicates if the node is selected, and is zero otherwise. More generally, we can define $f_{jk}(y)$, as the functional representing node $k$ in layer $j$, given as $N(j,k)$ in general. An active path link is defined between layer $j-1$ and layer $j$ for nodes $i$ and $k$ if node $N( j-1, i)$ and node $N(j, k)$ are both active. Thus, the marked path is derived from the set of selected nodes, and the collective measure must ensure that the path is feasible (\emph{i.e.} connected). A penalty term (P2) is introduced to account for path infeasibility, one stemming from invalid connections given the nodes selected. Note that no path calculation is necessary for a fully connected graph as each is automatically valid, and the penalty term P2 is zero by default. 


\vspace{3mm}

The node functional (for component $k$ in layer $j$) may also comprise a set of continuous variables X. Thus, the function $f_{jk}(X,y)$ can be stated for the general case for component $k$ in stage $j$. In this manner, each functional may represent a mixed-integer form, and the collective measure $F(\mathrm{X,Y})$ indicates the function over all nodes and stipulated variables. Similarly, a local constraint can be defined as $g_{jk}(X,y)$ at node $k$ in stage $j$, or as a global constraint $G(\mathrm{X,Y})$ over the entire decision horizon. In general terms, these constraints may be defined as linear, nonlinear, as equality or inequality, and specified by need. The definitions should be provisioned at the outset for each node $N(j, k)$ and for the overall objective along with constraints as required. These are anticipated definitions in the problem statement. Details can be found in the Algorithm section below. Note that path-dependent constraints can also be defined, as a function of the selected nodes. For example, to address variations arising from selected edges as per the intended problem specification. The evaluation of path infeasibility (P2) is one such example. Note that all constraints must be managed implicitly as penalties as the proposed binary solver is unconstrained. 

\vspace{3mm}

The formulated problem is posed as a binary nonlinear system. The binary variables stem from the selection of components in each stage by index (an integer) represented as a binary string. For example, with $M=2$ choices, a single bit gives two options [0 or 1], where 0 represents choice of component 1, while 1 indicates choice of component 2. For $M=4$, a two-bit string yields four choices [00, 01, 10, 11] with integer values [0, 1, 2, 3] by index. Here, index 0 is for component 1, index 1 is for component 2, and so on. The same notion applies for arbitrary size $M$, where the bit string size is first established as a function of the anticipated number of choices $M_j$ in stage $j$. If the resulting bit string range ($y_j^{max}$), is greater than $M_j-1$ (with the offset), then the excess index values will represent dummy nodes in the graph. A feasible or valid choice is one where the selected index $y_j$ is less than $M_j$. Thus, we introduce a penalty term (P1) for feasible component selection that serves to ensure only permitted components are selected by the solver. 
An elegant scheme is to stipulate a feasibility term $c_{jk}$ for component $k$ in layer $j$ \emph{a priori}, where the value is 0 if feasible and sufficiently high if invalid. P1 is then simply the sum of $c_{jk}$ over all active components with $y_{jk}=1$. Evaluation of path penalty P2 is a little more involved as each selected edge must first be identified and then validated as being viable. Thus, problems with a fully connected graphs are faster to evaluate as this evaluation step is avoided. Ultimately, a feasible path is one that connects valid components stage-by-stage over valid links between the selected nodes, with penalty terms P1 and P2 as zero. For invalid paths, P1 and P2 will be non-zero. The penalty terms are established and added to the composite penalty function. 
Similarly, additional user specified constraints (both local and global) must be accounted for in the combined measure. 
Thus, the overall problem can be stated as follows:

\begin{eqnarray}
\label{eqn:p2}
min~V(Z) = F(Z) + P_1(Z) + P_2(Z) + P_L(Z) + P_G(Z)   \\
Z \in \mathbb{B}^n 	\nonumber 	\\
z_i \in \{ 0, 1 \}	 ~~ \forall i \in [1~n]	\nonumber 
\end{eqnarray}

\vspace{2mm}

\noindent where $V$ indicates the collective penalty function and $Z$ is the representative set of binary variables of size $n$.
$P_1$ and $P_2$ are the component and path penalty terms, respectively.
In addition, $P_L$ represents the set of local constraints: 

\begin{equation}
P_L = \sum_{i, j, k} \gamma y_{jk}(Z)~ \mathrm{max} \{0, g_{jk}^i(Z) \}^2 
\end{equation}

\vspace{2mm}

\noindent where $\gamma$ is a penalty multiplier, $y_{jk}$ asserts if component $k$ in stage $j$ is active, with $i$-th local constraint $g_{jk}^i$.
Laslty, $P_G$ represents the set of global constraints:

\begin{equation}
P_G = \sum_{i} \gamma \mathrm{max} \{0, G^i(Z) \}^2 .
\end{equation}

\vspace{5mm}

A feasible path through a graph with full connectivity is shown in Figure \ref{fig:FullLink}, as compared to one with partial connectivity in Figure \ref{fig:PartialLink}, where the goal is to maximize the
collective value ensuing from each decision. Hence, the objective value decreases from 227 to 158 when certain links are no longer available. 
The robust stochastic solver is described in the next section, with algorithm details to follow.

\clearpage\newpage

\clearpage\newpage
\rhead{\tiny BINARY SOLVER}
\small
\section*{\sffamily \large\color{darkblue} Binary Nonlinear Solver }

The robust binary nonlinear stochastic solver (RSS) is designed to manage binary unconstrained nonlinear problems. Here, we assume that all continuous variables of interest are suitably encoded, 
and constraints are handled by penalty reformulation~(\ref{eqn:p2}). The problem is then stated as follows:

\begin{eqnarray}
\label{eqn:p1}
min~F(Z) 						\\
Z \in \mathbb{B}^n 	\nonumber 	\\
z_i \in \{ 0, 1 \}	 ~~ \forall i \in [1~n]	\nonumber 
\end{eqnarray}

\noindent where $Z$ is the set of binary variables of size $n$ and $F$ is the objective function of interest. Note that no stipulations are placed on the nature of this function. 

\vspace{3mm}

Specifically, $F(Z)$ represents the penalty function derived from the multi-stage problem definition, as discussed in the last section, where $Z$ of size $s_{leng}$ is the 
representative binary string composed by concatenation of decision strings representing all variables in the stipulated problem. 
The RSS serves to find the argument set $Z$ that optimizes the collective measure $F(Z)$, as per (\ref{eqn:p1}).
Conversion of the binary string $Z$ yields the index set corresponding to the selected component at each stage, along with any related continuous variables, if defined.
These additional variables are automatically appended in the construction of string $Z$, with length $b_{leng}$ for the binary set and length $c_{leng}$ for the continuous set.
All the variables are suitably extracted and decoded, per the problem definition, prior to evaluation.

\vspace{3mm}

The stochastic solution procedure is based on the notion of establishing a probability distribution measure for each binary variable in the problem.
This is achieved by sampling each binary variable commencing from a uniform distribution that is iteratively updated based on the information gathered from a random
sample set. In this manner, the distribution array Xd of size $n$ is sought that will indicate the likely probability of a binary variable taking a value of 1 (or 0 otherwise).
It is the intent, based on learning, that sampling the updated Xd over time will yield $\hat{Z}$ that optimizes $F(Z)$. Indeed, the process is sample dependent and as such, cannot guarantee that the solution is optimal. However, the procedure is robust, able to handle nonlinearity and the noted problem complexity. Empirically, the solver is shown to be highly effective for problems of moderate size ($\sim$$n$=1000). The method employs probabilistic learning with reinforcement and may be construed as a variant of population-based incremental learning used to enhance genetic algorithm operations~\cite{PBIL}.

\vspace{3mm}

The RSS employs an iterative schema based on stochastic generation and evaluation of sample paths. 
Here a path refers not to the selection of components in each stage of the multi-decision problem, but rather the evolution of the candidate control vector in $n$-space.
That is, it traces the improving control set $Z$ on a point-by-point basis with respect to the objective, where a given sample represents a solution over the problem graph\footnote{Note that
evaluation over a partial horizon would lead to a Monte-Carlo tree-search scheme.}. 
See for example Figure \ref{fig:Paths}. Thus, an evolving \emph{sample} path should not be confused with a \emph{solution} path in the graph. 

\vspace{3mm}

Initially, a collection of random candidates are generated, using an uniformly initialized distribution array Xd, where each represents a point in $n$ binary-space.
A Monte-Carlo (MC) search is performed from each given starting point and is akin to random path evolution in $n$ dimensions per sample. 
Each generated path $j$ in the sample set S is evaluated and its corresponding objective $f_j=F(Z_j)$ value is stored.
Next, the distribution array Xd is updated using the best generated paths based on ranked selection (typically, using the top 20\% by objective measure). 
The array Xd asserts a probability value for element $z_i$ to have a value of one in order to optimize the key functional $F(Z)$. 
Hence, Xd can be construed as an evolving estimate of the PDF of $Z$ for the intended goal. 
Future samples are thus generated using the updated distribution vector Xd. 
The process repeats for a given number of iterations for each initial candidate, where the best solution is updated to reflect to the best sample identified. 
Note that the sampling method is akin to epsilon-greedy approach, in which a small chance for exploration is always retained~\cite{SUTTON}. 
At the end of this cycle, information from all candidates is collated and an informed distribution Xd is derived over the group. 
A new cycle is then initiated using the updated distribution Xd (that was initially uniform) over all candidates for each bit $z_i$ in the string of size $n$. 
These candidates will thus begin to converge (enforcing \emph{exploitation}), yet maintain a small probability for random search (to enable \emph{exploration}).
Demonstrative early and late stage distribution evolution plots (for the Knapsack problem) are given in Figures \ref{fig:EarlyDist} and \ref{fig:LateDist}, respectively~\cite{KNAP}. 

\vspace{3mm}

In summary, the RSS process employs Monte-Carlo search with probabilistic-informed learning with reinforcement. 
Each design candidate carries its own distribution Xd that evolves by iteration.
Shared information over all candidates results in a collective update with a stipulated learning rate. 
The method is robust and empirically works well for moderate dimensionality ($\sim$$n=1000$)\footnote{The prototype was developed in MatLab and will benefit from a higher-level implementation.}.
The scheme is applicable to larger problem sizes but at greater sample cost.
This is an issue that impacts most reinforcement learning methods and is a possible area for future research.
Two test cases are presented next and the algorithmic details of the solver are presented in the section to follow.

\clearpage\newpage
\rhead{\tiny TEST CASES}
\small
\section*{\sffamily \large\color{darkblue} Test Cases}

\subsection*{\sffamily \small\color{darkblue} Test Case 1}

\begin{figure}[ht]
    \centering
    \includegraphics[width=0.5\textwidth]{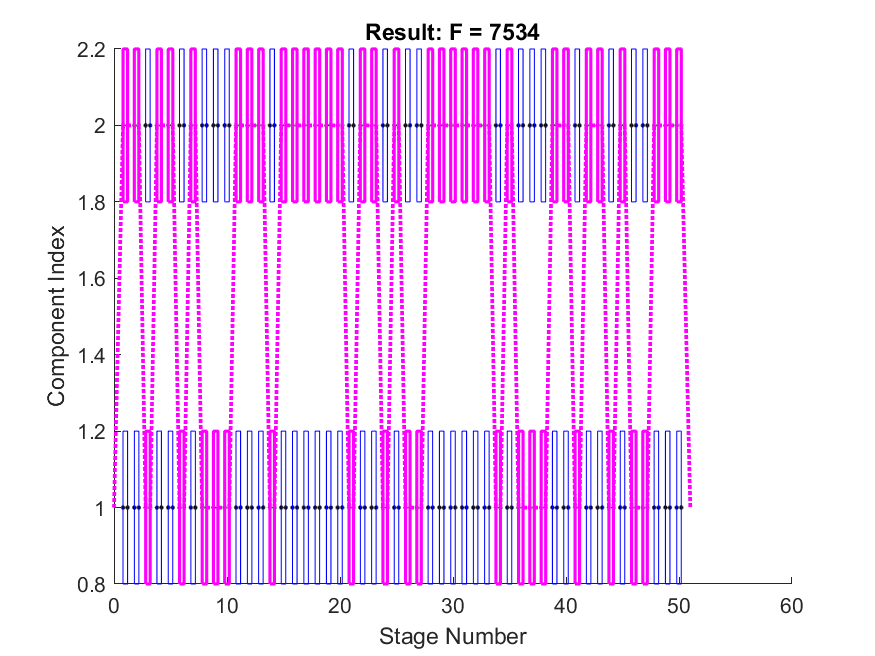}
    \caption{Solution path over 50 stages - Knapsack example~\cite{KNAP}. } 
	\footnotesize{\label{fig:Stages50} 50 stages with 2 components per stage. The binary system is $2^{50}$ with 1.1259e$^{15}$ choices.  } \\
    \label{fig:my_figure}
\end{figure}

\begin{figure}[ht]
    \centering
    \includegraphics[width=0.5\textwidth]{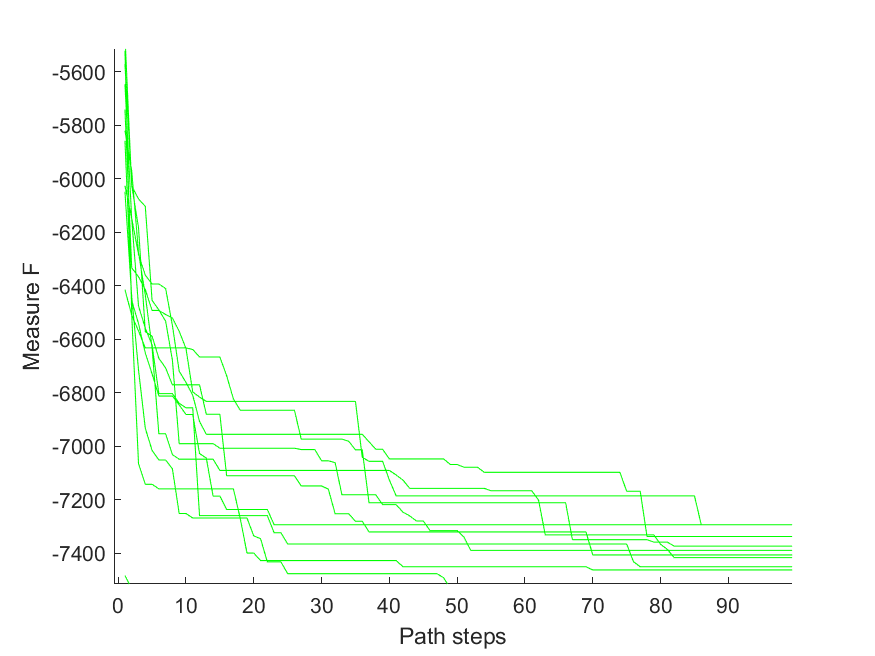}
    \caption{Candidate sample path evaluation - Knapsack example~\cite{KNAP}. } 
	\footnotesize{\label{fig:Paths} Shows improving sample paths for a set of candidates for a given iteration.  } \\
    \label{fig:my_figure}
\end{figure}

\begin{figure}[ht]
    \centering
    \includegraphics[width=0.6\textwidth]{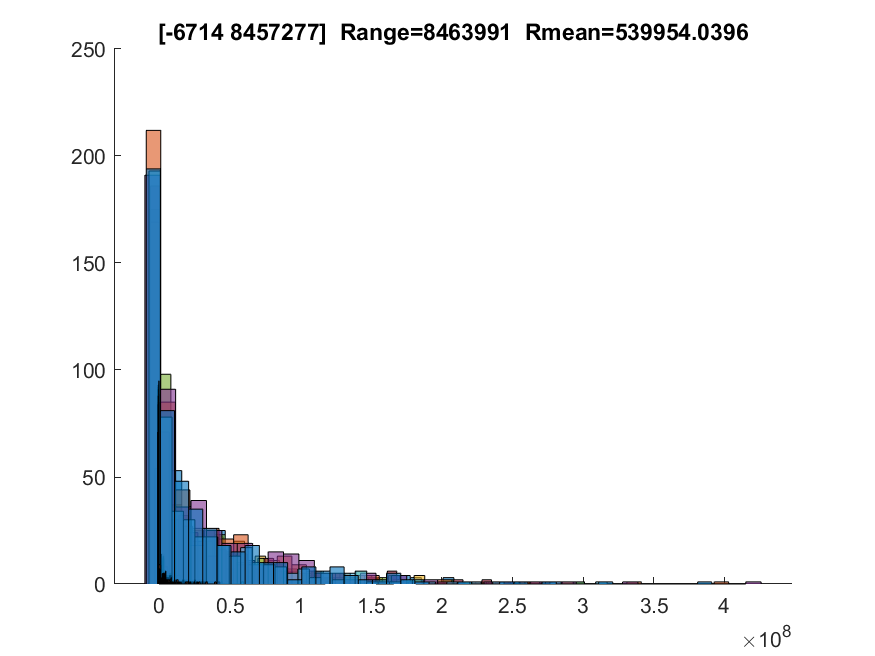}
    \caption{Early stage distribution estimation - Knapsack example~\cite{KNAP}. } 
	\footnotesize{\label{fig:EarlyDist} Shows distribution evolution (overlapping) with objective (x-axis) statistics: lower, upper, range and mean.  } \\
    \label{fig:my_figure}
\end{figure}

\begin{figure}[ht]
    \centering
    \includegraphics[width=0.6\textwidth]{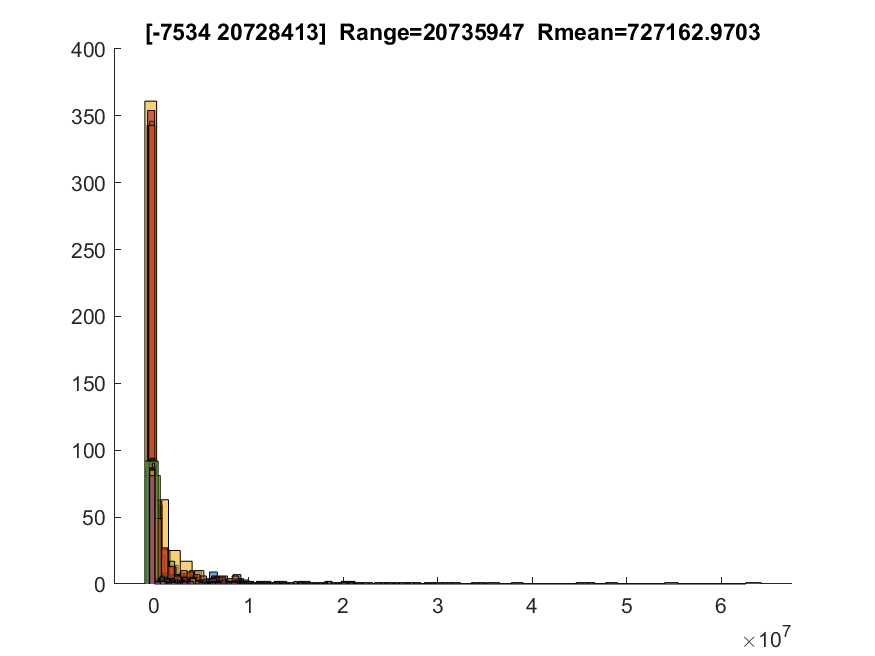}
    \caption{Late stage distribution estimation - Knapsack example~\cite{KNAP}. } 
	\footnotesize{\label{fig:LateDist} Shows distribution evolution (overlapping) with objective (x-axis) statistics: lower, upper, range and mean.  } \\
    \label{fig:my_figure}
\end{figure}

The proposed scheme is demonstrated using a classic Knapsack problem~\cite{KNAP}.
A decision-maker is given 50 items, each with a known utility and weight value, along with a weight capacity limit of 850.
The goal is to pick a number of items to maximize the utility subject to the weight limit.
This gives rise to $N=50$ stages, each with $M=2$ choices, to pick or not to pick an item.
The solution path, indicating an objective value $\hat F = 7534$, is shown in Figure \ref{fig:Stages50}.
Demonstrative candidate samples paths are shown in Figure \ref{fig:Paths}.
Lastly, early-stage and late stage evolving distributions are presented in Figures \ref{fig:EarlyDist} and \ref{fig:LateDist}, respectively.
Note that in Figure \ref{fig:LateDist} the lower bound is -7534 (the anticipated result when the solver is minimizing).

\clearpage\newpage
\subsection*{\sffamily \small\color{darkblue} Test Case 2}


In this example, we consider the $n$-dimensional Rosenbrock benchmark case~\cite{ROSEN}.
The problem is defined as follows:

\begin{eqnarray}
\mathrm{min}~f(X) = \sum_{i=1}^{n-1}   \left[  100 (x_{i+1} - x_i^2)^2 + (x_i - 1)^2  \right]		\\
							\nonumber  	\\
\mathrm{s.t.} -2 \le x_i \le 2			\nonumber		\\
		\forall i \in [1 \ldots n]		\nonumber		\\
		X \in \mathbb{R}^n 		\nonumber
\end{eqnarray}

\noindent where the known solution is given by: $\hat{X}$=[1 \ldots 1] with $\hat{F} = 0$.

\vspace{2mm}

This continuous problem can be posed as a multi-stage problem with $N$ stages and $M=1024$ choices in each.
The resulting 10-bit representation per stage gives a spatial granularity of 0.0039 (4/1024) over each dimension $x_i$ with bounds [-2~2].

The results are shown in Table \ref{ResTab} for increasing size $N$.
Clearly, this is a hard problem using binary encoding as simple variations result in dramatic changes in the continuous space impeding progress and convergence for larger instances\footnote{The use of Gray code may be more apt in this situation.}.
The method is intended for integer selectivity and is not ideal for continuous variable representation on the graph, \emph{i.e.}, with a large number of choices, $M=1024$.
Nonetheless, the example serves to demonstrate robustness under increasing model complexity.

\vspace{3mm}

\begin{table}[!h]
\begin{center}
\caption{Test Case 2 - Results}
\label{ResTab}
\begin{tabular}{|c|c|c|c|c|c|c|}
\hline
\small\emph{Stages}  	& \small\emph{Components}	& \small\emph{Graph Size}	&  \small\emph{Model} 		& \small\emph{Search}		& \small\emph{Solution} 	& \small\emph{Time}        \\
\small{$N$}  		& \small{$M$}			& \small{$N$x$M$}		&  \small\emph{Complexity} 	& \small\emph{Space}	& $F_{best}$ 		& \small\emph{mins}  	\\
\hline\hline
2		& 1024		&	 2048		&	 $2^{20}$	& 1.0486e$^6$		&	0.000096	&	1.71		\\
5		& 1024		&	 5120		&	 $2^{50}$	& 1.1259e$^{15}$		&	0.002684	&	2.98		\\
10		& 1024		&	10240 	&	 $2^{100}$	& 1.2677e$^{30}$		&	0.000867	&	4.81		\\
25		& 1024		&	25600	&	 $2^{250}$	& 1.8093e$^{75}$		&	0.037647	&	22.15		\\	
50		& 1024		&	51200	&	 $2^{500}$  & 3.2734e$^{150}$	&	16.82819     &	43.92		\\
\hline
\end{tabular}
\end{center}
\end{table}

\clearpage\newpage
\rhead{\tiny ALGORITHM DETAILS}
\small
\section*{\sffamily \large\color{darkblue} Algorithm Details }

The multi-stage decision framework and solver scheme are presented algorithmically in this section.

\vspace{3mm}

\noindent The high-level procedure, described in {\bf\texttt {Solution Procedure}}, comprises four steps. 
These entail defining the problem, setting algorithm parameters, parsing the problem and then solving the problem.

\vspace{3mm}

\noindent First, the problem is defined as per {\bf\texttt {Define MyProb}}. This includes stating the number of stages and the number of selectable components in each.
For non-fully connected graphs, the array of active links (ALIST) or the array of deactivated links (DLIST) must be provided. In addition, the number of continuous variables by stage and component are also defined  here by need. The problem definition is returned in the structure \texttt{PDEF}.

\vspace{3mm}

\noindent The algorithm parameters are returned in the structure \texttt{PARAS} as defined in {\bf\texttt {Set Paras}}.

\vspace{3mm}

\noindent The problem is then parsed in {\bf\texttt {Process MyProb}} to size the binary system and to automatically generate the problem graph.
The complete problem definition is returned in the structure \texttt{PROB}. 

\vspace{3mm}

\noindent Finally, the method {\bf\texttt {Run RSS}} calls the solver with the required input data and returns with the solution structure \texttt{SOLN}. 

\vspace{3mm}

\noindent The inner workings of the solver are relayed in the remaining procedures:
\newline\indent $\bullet$ {\bf\texttt {Run Cycle}} 
\newline\indent $\bullet$ {\bf\texttt {Evaluate Path}}
\newline\indent $\bullet$ {\bf\texttt {Sample Set}} 
\newline\indent $\bullet$ {\bf\texttt {Evaluate Set}} 
\newline\indent $\bullet$ {\bf\texttt {Evaluate Sample}} 
\newline\indent $\bullet$ {\bf\texttt {Iteration Evaluation}}
\vspace{1mm} 
\newline \noindent  The reader may review these methods qualitatively or in detail by interest. Note that the definitions are in line with those given in the preceding sections.

\vspace{3mm}

\noindent An example test problem (the Knapsack case) is presented in {\bf\texttt {Evaluate MyFunc}}.

\vspace{3mm}

\noindent Finally, note that for a given problem, the user need only furnish the definition and its evaluation method in {\bf\texttt {define MyProb}} and {\bf\texttt {evaluate MyFunc}}, respectively.

\vspace{8mm}

\noindent{\bf\texttt {Solution Procedure} }
\newline Set problem specific parameters
\newline Establish one-time pre-processed data
\newline Run problem instance:
\newline\indent Define the problem ({\bf\texttt {define~MyProb}}) $\to$ \texttt{PDEF} 
\newline\indent Set algorithm parameters ({\bf\texttt {set~Paras}}) $\to$ \texttt{PARAS}
\newline\indent Parse the problem ({\bf\texttt {process~MyProb}}) $\to$ \texttt{PROB}
\newline\indent Run the solver ({\bf\texttt {run~RSS}}) $\to$ \texttt{SOLN}
\newline\indent Show and plot results
\newline Create new instance

\clearpage\newpage

\noindent {\bf\texttt {Define MyProb} }
\newline Specify pertinent parameters: 
\newline\indent Number of decision stages, $N$ 
\newline\indent Maximization problem flag, $f_{max}$ ~[solver minimizes by default]
\newline\indent Selected function by index, $f_{eval}$ 
\newline\indent Penalty violation factor, $\gamma$
\newline\indent Actual solution, if known, $F_{act}$
\newline Establish model data given parameters
\newline Specify choices for each stage $j$ in $N$, \emph{e.g.}, \texttt{STAGE}($j$).choices = 2.
\newline Set array of active links (ALIST, size $n_{add}$) where flag connecting N($j-1$, $i$) to N($j$, $k$) is true.
\newline Set array of deactived links (DLIST, size $n_{del}$) where flag connecting N($j-1$, $i$) to N($j$, $k$) is false.
\newline Specify total number of continuous variables required ($f_{cvar}$).
\newline Set the continuous variable structure \texttt{CVAR} (if $f_{cvar} > 0$).
\newline \texttt{CVAR} indicates the continuous variable distribution by stage and component with desired ranges.
\newline Problem data is stored in the definition \texttt{PDEF}.
\newline Return: \texttt{PDEF}

\vspace{8mm}

\noindent {\bf\texttt {Set Paras} }
\newline Specify algorithm parameters: 
\newline\indent Number of cycles, $n_{cyc}$ = 6
\newline\indent Outer loop iterations, $k_{itns}$ = 12
\newline\indent Inner loop trial paths, $n_{path}$ = 48
\newline\indent Number of path steps, $p_{step}$ = 120
\newline\indent Initial sample size, $n_{s1}$ = 1e4 
\newline\indent Iteration sample size, $n_{s2}$ = 1e3
\newline\indent Relaxation factor, $\alpha$ = 0.01
\newline\indent Uniform sample flag, $f_{un}$ = 0
\newline\indent Fast update flag, $f_{upd}$ = 0
\newline\indent Plot distribution, $f_{plt1}$ = 0
\newline\indent Plot sample paths, $f_{plt2}$ = 0
\newline Return: \texttt{PARAS}

\clearpage\newpage

\noindent {\bf\texttt {Process MyProb} }
\newline \textit{Method parses specifed problem}
\newline function  \texttt{PROB} =  {\bf\texttt {process MyProb}}(\texttt{PDEF})
\newline Get parameters: $n_{add}$, $n_{del}$, $N$ and \texttt{STAGE}
\newline Initialize total bit count $M_{tot}$ and total number of combinations $M_c$.
\newline for each stage $k$ in $N$:
\newline\indent Get choices $m_k$
\newline\indent Establish number of bits required, $m_k^b$ given choices $m_k$.
\newline\indent Establish combinations $m_k^c$ over $m_k^b$ bits and the slack given $m_k$ required choices.
\newline\indent Set selected index range [1~$m_k$]
\newline\indent Update cumulative values: $M_{tot}$ and $M_c$
\newline\indent Store data in \texttt{STAGE}(k) 
\newline Establish the string block for each stage $k$ over the total length $M_{tot}$ by index [$k_{st}$~$k_{ed}$]
\newline Establish if the graph is fully connected and set flag IsCon. 
\newline Set problem structure comprising:
\newline\indent Updated \texttt{STAGE} information
\newline\indent Combinations by stage, $m_k^c$
\newline\indent Binary system length, $b_{leng}$ = $M_{tot}$
\newline\indent Binary problem flag, IsBin
\newline\indent Fully connected flag, IsCon
\newline\indent Generated problem graph \texttt{GRPH} given stage information and connectivity map.
\newline\indent Variable range definition with LU arrays.
\newline\indent Continuous variables and ranges, with cumulative string size $c_{leng}$.
\newline\indent Total string length, $s_{leng}$ = $b_{leng}$ + $c_{leng}$.
\newline\indent Full problem combinations, $M_{tot}$ given size $s_{leng}$. 
\newline Return: \texttt{PROB}

\vspace{7mm}

\noindent {\bf\texttt {Run RSS} }
\newline \textit{Execute the solver given the parameter and problem statement}
\newline function  \texttt{SOLN} =  {\bf\texttt {run RSS}}(\texttt{PARA}, \texttt{PROB})
\newline Get parameters: $n_{cyc}$, $\alpha$ and $f_{eval}$
\newline Initialize empty arrays Xin, Vin, Vbest, PV and PF. Set Fin=0, etime=0, Fbest=1e6.
\newline for each cycle $j$ in $n_{cyc}$:
\newline\indent [Xd, V, F, et] = {\bf\texttt {run cycle}}(\texttt{PARAS}, \texttt{PROB}, Xin, Vin, Fin)
\newline\indent Store results PV [$N$ $m$], PF [$m$ 1] and update etot, with stages $N$ and samples $m$.
\newline\indent if(F $<$ Fbest):  update Fbest=F and Vbest=V
\newline end
\newline Perform final evaluation step over all data in PV:
\newline Xd =  {\bf\texttt {iteration evaluation}}(PV, $\alpha$/10)
\newline [Xsol, Vsol, Fsol, et] =  {\bf\texttt {run cycle} }(\texttt{PARAS}, \texttt{PROB}, Xd, Vbest, Fbest)
\newline Update elapsed time, etot
\newline Set decoded solution D [1 $N$] with path link data PLINK
\newline Evaluate at solution:  RESP = {\bf\texttt {eval set}}(Vsol, \texttt{PARAS}, \texttt{PROB})
\newline Get Ysol, X, XS, PLINK from RESP
\newline Store solution: PF, Xsol, Vsol, Fsol, Ysol, Csol, X, XS, PLINK, etot, optgap
\newline Return: \texttt{SOLN}

\clearpage\newpage

\noindent {\bf\texttt {Run Cycle} }
\newline \textit{Evaluation call given input parameters (with existing distribution array Xd) }
\newline function [Xd, Vbest, Fbest, et] =  {\bf\texttt {run cycle}}(\texttt{PARAS}, \texttt{PROB}, Xdin, Vin, Fin)
\newline Get parameters: 
\newline\indent Relaxation factor, $\alpha$
\newline\indent Outer loop iterations, $k_{itns}$ 
\newline\indent Inner loop trial paths, $n_{path}$ 
\newline\indent Problem dimensionality, $s_{leng}$ 
\newline\indent Maximization problem flag, $f_{max}$
\newline\indent Fast update flag, $f_{upd}$
\newline Get length of Xdin as $n_{test}$
\newline if($n_{test}$ = 0):  initialize Xd as uniform distribution [$N$ 1], Vbest=[~] and Fbest=1e6
\newline if($n_{test}$ $>$ 0):  set Xd=Xdin [$N$ 1], Vbest=Vin [$N$ 1] and Fbest=Fin
\newline if($f_{max}$):  fmult=-1 
\newline for $j$ in $k_{itns}$:
\newline\indent set PTH struct [~]
\newline\indent for all paths $k$ in $n_{path}$:
\newline\indent\indent [v, f] =  {\bf\texttt {evaluate path}}(\texttt{PARAS}, \texttt{PROB}, Xd)
\newline\indent\indent Store v and f in PTH($k$)
\newline\indent end
\newline\indent Initialize array V=[~]
\newline\indent for each path $k$ in $n_{path}$:
\newline\indent\indent Get v and f in PTH($k$)
\newline\indent\indent if(f $<$ Fbest):  Fbest = f  and Vbest = v
\newline\indent\indent Store v in V [$N$ $n_{path}$]
\newline\indent end
\newline\indent Update Xd for next iteration $j$:  Xd = {\bf\texttt {iteration evaluation}}(V, $\alpha$, Xd)
\newline\indent if(fupd): fast update: Xd = {\bf\texttt {iteration evaluation}}(Vbest, $\alpha$, Xd)
\newline end 
\newline Set final Xd = {\bf\texttt {iteration evaluation}}(Vbest, $\alpha$, Xd)
\newline Return: Xd, Vbest, Fbest

\clearpage\newpage

\noindent {\bf\texttt {Evaluate Path} }
\newline \textit{Generate and evaluate a path with given distribution }
\newline function [Ybest, Fbest] =  {\bf\texttt {eval path} }(\texttt{PARAS}, \texttt{PROB}, Xd)
\newline Get parameters: 
\newline\indent Initial sample size, $n_{s1}$
\newline\indent Iteration sample size, $n_{s2}$ 
\newline\indent Relaxation factor, $\alpha$
\newline\indent Number of path steps, $p_{step}$
\newline\indent Uniform sample flag, $f_{un}$
\newline\indent Initialize STORE = [~]
\newline for $k$ in $p_{step}$:
\newline\indent Set X = Xd
\newline\indent if($k=1$): S = {\bf\texttt {sample set}}(X, $n_{s1}$, $f_{un}$)
\newline\indent if($k>1$): S = {\bf\texttt {sample set}}(X, $n_{s2}$, 0)
\newline\indent [F, Y, RESP] =  {\bf\texttt {eval set}}(S, \texttt{PARAS}, \texttt{PROB})
\newline\indent if($k=1$):  initialize Fbest=F and Ybest=Y
\newline\indent if(F $<$ Fbest):  Fbest=F and Ybest=Y
\newline\indent Update Xd = {\bf\texttt {iteration eval}}(Y, $\alpha$, X) given pop Y
\newline\indent Update STORE:  $k$, $\alpha$, $F$, $F_{mean}$, $F_{lwr}$, $F_{upr}$
\newline end
\newline Plot generated path (steps vs measure)
\newline Return: Ybest, Fbest

\vspace{1cm}

\noindent {\bf\texttt {Sample Set} }
\newline \textit{Generate m sample paths S [m n] given distribution array X and the uniform sampling flag }
\newline function S =  {\bf\texttt {sample set}}(X, $m$, $f_{un}$)
\newline Get dimensions, $n$ 
\newline Initialize random variable array, RV [$m$ $n$]
\newline Initialize output sample array S [$m$ $n$] as ones
\newline for each dim $k$ in $n$:
\newline\indent $x$ = X($k$)
\newline\indent if($f_{un}$): $x$ = 0.5
\newline\indent for each sample $j$ in $m$:
\newline\indent\indent $\beta$ = RV($j$,$k$)
\newline\indent\indent if($\beta < x$):  S($j$, $k$) = 0
\newline\indent\indent if($j=1$) and ($x \ge$ 0.5):  S($j$, $k$) = 0
\newline\indent end
\newline end
\newline Return: sample set S

\clearpage\newpage

\noindent {\bf\texttt {Evaluate Set} }
\newline \textit{Evaluates the function value F for each sample in B in BPOP and returns the best }
\newline function [Fbest, Ybest, RESP] =  {\bf\texttt {evaluate set}}(BPOP, \texttt{PARAS}, \texttt{PROB})
\newline Get number of samples, $m$ in BPOP [$m$ $n$]
\newline Initialize the solution array, FS [$m$ 1]
\newline for each sample $j$ in $m$:
\newline\indent Get $j$-th sample B [1 $n$]
\newline\indent Evaluate sample FS($j$) =  {\bf\texttt {evaluate sample}}(B, \texttt{PROB})
\newline end
\newline Get population statistics (including Fmean)
\newline Establish Fbest and corresponding Bbest
\newline Get properties for Ybest:  [Fchk, Y, X] =  {\bf\texttt {evaluate sample}}(Bbest, \texttt{PROB})
\newline Establish path link information PLINK given design D for Ybest
\newline Store solution responses RESP:  Fmean, Y, X, PLINK
\newline Return: Fbest, Bbest, RESP

\vspace{1cm}

\noindent {\bf\texttt {Evaluate Sample} }
\newline \textit{Evaluate binary sample B giving objective F, integer set Y, continuous set X }
\newline function [F, Y, X] =  {\bf\texttt {evaluate sample}}(B, \texttt{PROB})
\newline Get parameters: 
\newline\indent Penalty factor, $\gamma$
\newline\indent Binary problem flag, IsBin
\newline\indent Fully connected problem flag, IsCon
\newline\indent Binary problem size, $b_{leng}$
\newline\indent Continuous variable size, $f_{cvar}$
\newline\indent Maximization flag, $f_{max}$
\newline\indent Stage data structure, \texttt{STAGE}
\newline Extract binary set BY from B of size $b_{leng}$
\newline if(IsBin=1): integer selection directly, Y = BY
\newline if(IsBin=0): integer selection by string conversion, Y = {\bf\texttt {bin2control}}(BY, \texttt{STAGE})
\newline Initialize continuous variable set X = [~]
\newline if($f_{cvar}>0$):  X = {\bf\texttt {bin2cont}}(B, \texttt{PROB}) - convert appended string (size $c_{leng}$) for continuous variables.
\newline Initialize the objective and penalty terms: f = P1 = P2 = 0
\newline if(IsBin=0): get component penalty P1 = {\bf\texttt {eval cmpt penalty}}(Y, \texttt{STAGE}, 10$\gamma$)
\newline if(IsCon=0): get path penalty P2 =  {\bf\texttt {eval path penalty}}(Y, \texttt{PROB})
\newline if(P1=0) and (P2=0): f = {\bf\texttt {evaluate myfunc}}(X, Y, \texttt{PROB})
\newline if($f_{max}$): set f = -f
\newline Set penalty function F = f + P1 + P2
\newline Return: F, Y, X

\clearpage\newpage

\noindent {\bf\texttt {Iteration Evaluation} }
\newline \textit{Establishes distribution over the set of m samples in V [n m] for a given relaxation and learning rate. }
\newline function Xd =  {\bf\texttt {iteration eval}}(V, $\alpha$, Xold)
\newline Get parameters: 
\newline\indent Relaxation factor, $\alpha$
\newline\indent Learning rate, $\lambda$ = 0.8 
\newline Get dimensions, $n$
\newline Get samples, $m$
\newline Estimate the row mean for each dimension $n$ over all $m$ samples giving S [$n$ 1]
\newline Initialize the distribution array Xnew [$n$ 1] as zero
\newline for each dimension $k$ in $n$:
\newline\indent Get estimate, $s_k$ = S($k$)
\newline\indent Set probability measure Xnew($k$) = 1 - $s_k$
\newline\indent if($s_k$=0): Xnew($k$) = 1-$\alpha$  (relaxed 0, with tiny probability for 1)
\newline\indent if($s_k$=1): Xnew($k$) = $\alpha$ (relaxed 1, with tiny probability for 0)
\newline end
\newline Set distribution with learning rate, Xd = Xold + $\lambda$(Xnew - Xold)
\newline Return: Xdt

\vspace{1cm}

\noindent {\bf\texttt {Evaluate MyFunc} }
\newline \textit{Knapsack example with binary input B [1 n] and problem struct PROB}
\newline function f =  {\bf\texttt {eval myfunc}}(B, \texttt{PROB})
\newline Get parameters: 
\newline\indent Utility array, U [$n$ 1] 
\newline\indent Weight array, W [$n$ 1]
\newline\indent Capacity limit, C
\newline\indent Penalty multiplier, $\gamma$
\newline Get utility sum, $u$ = B~U
\newline Get weight sum, $w$ = B~W
\newline Get constraint term, g = max(0, $w$-C)
\newline Set objective measure, f = u -  $\gamma~g^2$  [to be maximized]
\newline Return: f

\clearpage\newpage
\rhead{\tiny SUMMARY}
\small
\section*{\sffamily \large\color{darkblue} Summary}

The key points of interest are as follows:

\vspace{3mm}

\noindent $\bullet$ A framework for multi-stage decision problems represented on an underlying graph. 
\newline\noindent $\bullet$ A robust nonlinear binary solver based on incremental probabilistic learning with reinforcement to treat the combinatorial definition.
\newline\noindent $\bullet$ Method is applicable to various cases and allows treatment of complex graph arrangements by need.
\newline\noindent $\bullet$ Allows partial connectivity between stages and components, with varying number of choices by stage.
\newline\noindent $\bullet$ Deals with nonlinear relations between components and permits consideration of alternative pathways.
\newline\noindent $\bullet$ Problems are easier to define as the method does not require linearization or the use of MINLP methods with disjunctions to manage options.
\newline\noindent $\bullet$ The solution path is an abstraction for evaluation. The decision order is fungible. 
\newline\noindent $\bullet$ Problem size dictates computational cost. A large problem requires more samples for exploration. 
\newline\noindent $\bullet$ Solution is not  guaranteed to be optimal. It is sample and paramterization dependent.
\newline\noindent $\bullet$ The prototype will benefit from implementation in a higher language.
\newline\noindent $\bullet$ Testing and validation on other suitable problems is desirable to demonstrate utility\footnote{This method was successfully applied to an Energy Harvesting problem. The results are not reported here.}.

\clearpage\newpage
\rhead{\tiny REFERENCES}
\bibliographystyle{unsrt}
\small


\vspace{50mm}
\section*{Acknowledgement}
This work was undertaken in 2025 at Schlumberger-Doll Research (SDR), One Hampshire Street, Cambridge, MA 02139.

\end{document}